\documentclass[12pt,a4paper]{article}
\usepackage{amssymb}

\newcommand{\beq}{\begin{equation}}
\newcommand{\eeq}{\end{equation}}
\newcommand{\bb}{\ifmmode 
  \else \leavevmode\unskip\penalty9999 \hbox{}\nobreak\hfill
  \fi
  \quad\hbox{\hskip.5em\vrule width.5em height.6em depth.05em\hskip.1em}
  \vspace{2mm}}

\newtheorem {Lem} {Lemma}

\def \rit {\mathbb{R}} 

\def \nit {\mathbb{N}}

\title{A remark on compact $H$-surfaces into $\rit^3$}
\author{Yuxin Ge and Fr\'ed\'eric H\'elein}
\date{}

\begin{document}
\maketitle
\def\cqfd{%
\mbox{ }%
\nolinebreak%
\hfill%
\rule{2mm} {2mm}%
\medbreak%
\par%
}
\def \vs{\vspace*{0.2cm}}
\def \hs{\hspace*{0.6cm}}
\def \ds{\displaystyle}

\def \E{\frac {{\| \nabla \varphi \|}_2^2} {{{\| \nabla a \|}_2^2}{{\| \nabla b \|}_2^2}}}

\section{Introduction}Let $\Omega$ be a smooth and bounded domain in $\rit^2$. We
consider the following system
\begin{equation}
\triangle u = u_x \wedge u_y, \hs \mbox{in } \Omega,
\label{eq1}
\end{equation}
where $u\in C^2(\Omega;\rit^3)$ and subscripts denote partial differentiation 
with respect to coordinates. This equation characterizes surfaces of constant mean
curvature $H=\frac{1}{2}$ in $\rit^3$ in conformal representation. More
precisely, any non constant smooth map $u$ which is a solution of (\ref{eq1}) and of
the conformality condition
\begin{equation}
\omega:=(|u_x|^2- |u_y|^2- 2i\langle u_x, u_y\rangle )dz\otimes dz =0,
\hs
\mbox{in }
\Omega,
\label{eq1bis}
\end{equation}
(here $dz =dx+idy$) parametrizes a branched immersed constant mean curvature surface
in $\rit^3$. For that reason (\ref{eq1}) is called the $H$-system.
The complex
tensor $\omega$ which appears in (\ref{eq1bis}) is called the Hopf differential (see
\cite{J}).  The first existence result for
solutions to (\ref{eq1}) and (\ref{eq1bis}) was proved by H. Wente in \cite{H1}.

In \cite{He2}, the second author proposed a new variational approach for
finding a solution to (\ref{eq1}). For any pair of functions $a,b\in H^1(\Omega)$, we denote by
$\varphi:=\widetilde{ab}$ the unique solution in
$H^1(\Omega)$ of the Dirichlet problem
\begin{eqnarray}
\left\{
\label{eq2}
\begin{array}{rll}
- \triangle \varphi & =\{a,b\}, & \mbox{in } \Omega\\
\varphi & = 0, & \mbox{on } \partial \Omega,
\end{array}
\right.
\end{eqnarray}
where $\{a,b\}= {a_x}{b_y} - {a_y}{b_x}$. By a result of Brezis and Coron \cite{BC}
based on an idea due to H. Wente \cite{H1} \cite{W2}, we know that $\varphi$ is
continuous on $\bar\Omega$ and
\begin{eqnarray}
\label{eq3}
\| \varphi \|_{L^\infty(\Omega)} + \| \nabla\varphi \|_{L^2(\Omega)}\leq C_0(\Omega)
\| \nabla a \|_{L^2(\Omega)}\| \nabla b \|_{L^2(\Omega)}.
\end{eqnarray}
Thus the following energy functional makes sense
$$
E(a,b,\Omega) = \frac{\| \nabla a \|^2_{L^2(\Omega)}+\| 
\nabla b \|^2_{L^2(\Omega)}}{2\| \nabla \varphi \|_{L^2(\Omega)}},
\mbox{ defined for
}a,b\in H^1(\Omega)\setminus \{0\}.
$$
The Euler-Lagrange equation satisfied by the critical points of this functional was
derived in \cite{He2}.  Through the substitution $u:=(\lambda a, \lambda b, \lambda^2
\varphi)  $ for
$\lambda =  -\sqrt{\frac{{\| \nabla a \|}^2_2+{\| \nabla b \|}^2_2}{2{\| \nabla
\varphi \|}^2_2}}$, this equation coincides with (\ref{eq1}). The boundary
conditions are 
\begin{equation}
\varphi=\frac{\partial a}{\partial n}=\frac{\partial b}{\partial n} =0 \hs \mbox{on } \partial\Omega,
\label{eq4}
\end{equation}
where $n=(n_1,n_2)$ is the normal vector on $\partial\Omega$. Moreover,
 in \cite{Ge1} the first author showed that the Hopf differential $\omega$ is
holomorphic and satisfies the boundary condition $Im(\omega \nu^2)=0$, where
$\nu=n_1+in_2$. This implies in particular that $\omega$ vanishes - i.e.
(\ref{eq1bis}) is true - if $\Omega$ is simply connected.\\

An important property of our problem is that the functional $E$ and its critical points are preserved by conformal
transformations of the domain $\Omega$ (see \cite{He2}). So this variational problem depends only on
the complex structure of $\Omega$ and hence it also makes sense  to consider the problem on a Riemann surface. The 
boundary conditions (\ref{eq4}) allow us to construct a solution of (\ref{eq1}) from
a compact oriented Riemannian surface into
$\rit^3$ by gluing together two copies of
$\Omega$. More precisely, we construct $N:=\Omega\cup_{\partial\Omega}\tilde\Omega$,
where $\tilde\Omega$ is a copy of $\Omega$, provided with opposite orientation and
define a $C^\infty$ map
$\tilde u:N\rightarrow \rit^3$ by $\tilde u=u$ on $\Omega$ and $\tilde
u=(\lambda a,
\lambda b, -\lambda^2 \varphi)$ on $\tilde\Omega$. This map is a solution of the
$H$-system (\ref{eq1}) and its Hopf differential is holomorphic. Were this
differential to vanish, i.e. if $\omega =0$, then we would obtain a constant mean
curvature branched immersion. Recall that Wente \cite{Wente} (see also
\cite{Ab})  constructed an immersed constant mean curvature torus, which enjoys
an invariance under an orthogonal symmetry with respect to a plane. Thus it has the
form $\Omega\cup_{\partial\Omega}\tilde\Omega$ as above, where $\Omega $ is some
annulus. This motivates the search for critical points of $E$. In \cite{Ge1}, an
existence result was derived for a perforated domain, provided the holes
are small enough (in the same spirit as in \cite{C}). Here we address the problem of a
one-connected domain $\Omega$, i.e. a domain of the form $U\setminus \bar V$ where
$U$ and $V$ are smooth bounded simply connected open sets and $\bar V\subset U$,
without smallness assumption on the hole $V$. Because any such domain is
conformally equivalent to a radially symmetric annulus \cite{A} and thanks to
the invariance of the variational problem under conformal transformations, we shall restrict
 ourself to annuli without loss of generality. Our method relies on a minimization 
procedure on a subset of $H^1\times H^1$, which is equivariant with respect to some finite
group.\\
 
{\bf Main Theorem.} {\it Consider the annulus $\Omega:=\{(x,y)\in \rit^2,\hs
r_0<r=\sqrt{x^2+y^2} < 1\}$ with $0<r_0<1$. Then, there exists
a critical point $(a,b)\in H^1(\Omega, \rit^2)\setminus \{(0,0)\}$ of the energy
functional $E$. Moreover $(a,b)$ and $\varphi:=\widetilde{ab}$ are smooth and 
 satisfy the boundary conditions (\ref{eq4}). Thus there exists a real number 
$\lambda\neq 0$ such that the map $u:=(\lambda a,\lambda b, \lambda^2\varphi)$
is a solution of (\ref{eq1}). Lastly the Hopf differential $\omega$ of $u$ has
the form $\omega={\tau\over z^2}dz\otimes dz$ for some real number $\tau$.}\cqfd

Unfortunately we are not able to prove that the map $u$ is conformal. We expect that it should
be so for some values of $r_0$. Indeed the parameter $\tau$ characterizing $\omega$
should vary with $r_0$ (since the set of holomorphic quadratic differentials is,
roughly speaking, the dual space of Teichm\"uller space). But we are still far from
understanding how $\tau$ could be related to $r_0$.\\

\section{The Euler-Lagrange Equation}
First, we note that $\Omega$ is invariant under rotations. We define 
$F_m=\{\Theta=(a,b)\in H^1(\Omega)\times H^1(\Omega),\hs \Theta\circ A = A\circ
\Theta\}$ and we will prove existence of a minimum of $E$ for $m$ large enough, where
$A$ is the rotation of angle $\frac{2\pi}{m}$ in $\rit^2$.\\

\begin{Lem} 
\label{lemma1}
Assume that $(a,b)\in F_m$. Then, the unique solution $\varphi$ of
(\ref{eq2}) is invariant under $A$, that is,
$$
\varphi\circ A = \varphi.
$$\cqfd
\end{Lem}
{\em Proof.} Clearly, we have
$$
d\Theta = A^{-1}\cdot (d\Theta)\circ A \cdot A.
$$
Thus,
$$
\{a,b\}=det(d\Theta) = det (A^{-1})\cdot[det (d\Theta)]\circ A \cdot det(A) = [det (d\Theta)]\circ A = \{a,b\}\circ A,
$$
since $det(A)=det(A^{-1})=1$. On the other hand, the unique solution $\varphi$ of (\ref{eq2}) is also the
unique minimum of the following energy functional $E_1$:
$$
E_1(\psi) = \frac{1}{2}\int_\Omega |\nabla \psi|^2- \int_\Omega\{a,b\}\psi,\hs\mbox{defined for all } \psi\in H^1_0(\Omega).
$$
Obviously, $E_1(\varphi\circ A) = E_1(\varphi)$. By the uniqueness of the minimizer,
we deduce that 
$$
\varphi\circ A = \varphi.
$$
\cqfd

In order to get the Euler-Lagrange equation of $E$, we first recall  a technical
lemma inspired by the work in \cite{H1} and proved in the Appendix in \cite{BC}.\\

\begin{Lem} 
\label{lemma2}
If $\varphi \in H^1(\Omega) \cap
L^\infty(\Omega)$, $a \in H^1(\Omega) \cap L^\infty(\Omega)$, $b \in H^1(\Omega)$
and $\varphi a = 0$ on $\partial\Omega$, then we have
$$
\int_{\Omega}\varphi\{a, b\} = \int_{\Omega}a \{b, \varphi\}.
$$
\cqfd
\end{Lem}

The following result shows that critical points of $E$ on $F_m$ are also critical
points of $E$ on $H^1\times H^1$.\\

\begin{Lem} 
\label{lemma3}
Assume that $H=(a, b) \in F_m$ is a minimizer of $E$ on $F_m$. Then\\
1) there exists $\lambda \in \rit^*$ such that $\Psi = (a_1, b_1, \varphi_1) = (\lambda a, \lambda b, \lambda^2\varphi)$ satisfies equation (\ref{eq1}).\\
2) $\Psi$ verifies the boundary conditions (\ref{eq4}).\\
3) $\ds{\int_\Omega \nabla a \cdot\nabla b = 0}$.\\
4) $\|\nabla a \|_{L^2} = \|\nabla b \|_{L^2}$.\\
5) there exists $c \in \rit$ such that
$$
\ds{\langle \partial_z \Psi, \partial_z \Psi \rangle = \frac{c}{z^2}},
$$ 
where $\partial_z  =\ds{ \frac{1}{2}(\partial_x - i \partial_y)}$.\\
6) $\Psi$ is regular on $\bar\Omega$.
\cqfd
\end{Lem}
{\em Proof.} The proof is very similar to the proof of Theorem 3.2 in \cite{Ge1}.
We just need to adapt it to our equivariant setting.\\
Let
$\Lambda=(\alpha,\beta)\in F_m$. Denote by
$\psi
$ the unique solution of the following equation
\begin{eqnarray}
\left\{
\begin{array}{rll}
- \triangle \psi & =\{\alpha,b\} + \{a,\beta \}, & \mbox{in } \Omega\\
\psi & = 0, & \mbox{on } \partial \Omega.
\end{array}
\right.
\end{eqnarray}
We claim that $\psi\circ A = \psi$. Indeed, note that,
$$
\{\alpha,b\}+\{a,\beta\}= det (dH + d\Lambda) -det (dH) -det (d \Lambda)
$$
and thus as in Lemma 1, we deduce that
\begin{eqnarray}
\label{eq711}
(\{\alpha,b\} + \{a,\beta \})\circ A = \{\alpha,b\} + \{a,\beta \}.
\end{eqnarray}
Hence, by the same argument as before, we establish the claim. Now set $\Theta_t =
\Theta + t \Lambda$. Clearly, $\Theta_t\in F_m$. A direct calculation leads to
\begin{eqnarray}
\label{eq8}
\begin{array}{lll}
E(a_t, b_t, \Omega) &=& \frac{\ds{{\| \nabla a \|}^2_2+{\| \nabla b \|}^2_2+2t\int_{\Omega}({ \nabla a }\cdot{ \nabla \alpha }+ { \nabla b }\cdot{ \nabla \beta })+O(t^2)}}{\ds{2\sqrt{{\| \nabla \varphi \|}^2_2+2t\int_{\Omega}{ \nabla \varphi }\cdot{ \nabla \psi }+O(t^2)}}}\\
&=& \frac{\ds{{\| \nabla a \|}^2_2+{\| \nabla b \|}^2_2+2t\int_{\Omega}({ \nabla a }\cdot{ \nabla \alpha }+ { \nabla b }\cdot{ \nabla \beta })+O(t^2)}}{\ds{2\left({\| \nabla \varphi \|}_2+\frac{t}{{\| \nabla \varphi \|}_2}\int_{\Omega}{ \varphi (\{\alpha,b\} + \{a,\beta \})}+O(t^2)\right)}}\\
&=& E(a, b, \Omega)\left(1- \ds{\frac{t}{{\| \nabla \varphi \|}^2_2}}\ds{\int_{\Omega}{ \varphi (\{\alpha,b\} + \{a,\beta \})}}\right.\\
&&\left.
 +\ds{ \frac{2t}{{\| \nabla a \|}_2^2 + {\| \nabla b \|}_2^2}\int_{\Omega}({ \nabla a }\cdot{ \nabla \alpha }+ { \nabla b }\cdot{ \nabla \beta })+ O(t^2)}\right).
\end{array}
\end{eqnarray}
\noindent For fixed $\theta_0\in[0,2\pi]$, we consider the domain
$\Omega_{\theta_0} =\{(x,y),\hs r_0<r<1\mbox{ and }\theta_0<
\theta<\theta_0+\frac{2\pi}{m}\}$. It follows from (\ref{eq711}) and Lemma 1 that
$$
[\varphi(\{\alpha,b\} + \{a,\beta \})]\circ A = \varphi(\{\alpha,b\} + \{a,\beta \}),
$$
which implies
\begin{eqnarray}
\label{eq9}
\int_\Omega \varphi(\{\alpha,b\} + \{a,\beta \}) = m \int_{\Omega_{\theta_0}} \varphi(\{\alpha,b\} + \{a,\beta \}).
\end{eqnarray}
On the other hand, we have
\begin{eqnarray}
\begin{array}{ll}
({ \nabla a }\cdot{ \nabla \alpha }+ { \nabla b }\cdot{ \nabla \beta })\circ A&=
tr([d\Theta]\circ A\cdot [d\Lambda]^t\circ A)\\
&= tr(A\cdot d\Theta\cdot A^t\cdot A\cdot d\Lambda^t\cdot A^t)\\
&= tr(A\cdot d\Theta\cdot  d\Lambda^t\cdot A^t)\\
&= tr(d\Theta\cdot  d\Lambda^t\cdot A^t\cdot A)\\
&= { \nabla a }\cdot{ \nabla \alpha }+ { \nabla b }\cdot{ \nabla \beta },
\end{array}
\end{eqnarray}
that is,
\begin{eqnarray}
\label{eq10}
\int_{\Omega}({ \nabla a }\cdot{ \nabla \alpha }+{ \nabla b }\cdot{ \nabla \beta })=m\int_{\Omega_{\theta_0}}({ \nabla a }\cdot{ \nabla \alpha }+{ \nabla b }\cdot{ \nabla \beta }).
\end{eqnarray}
Combining (\ref{eq8}) to (\ref{eq10}), we obtain 
$$
\int_{\Omega_{\theta_0}}({ \nabla a }\cdot{ \nabla \alpha }+{ \nabla b }\cdot{ \nabla \beta })=\frac{{\| \nabla a \|}^2_2+{\| \nabla b \|}^2_2}{2{\| \nabla \varphi \|}^2_2}\int_{\Omega_{\theta_0}}\varphi (\{\alpha,b\}+\{a,\beta\}).
$$
In particular, if we set $\alpha$, $\beta \in C^\infty_0(\Omega_{\theta_0})$, we deduce from Lemma 2 that
\begin{eqnarray}
\left\{
\begin{array}{ll}
- \triangle  a =\vs\ds{ \frac{{\| \nabla a \|}^2_2+{\| \nabla b \|}^2_2}{2{\| \nabla \varphi \|}^2_2} \{ b, \varphi \}},&\mbox{in }\Omega_{\theta_0}\\
- \triangle  b =\ds{ \frac{{\| \nabla a \|}^2_2+{\| \nabla b \|}^2_2}{2{\| \nabla \varphi \|}^2_2}  \{ \varphi, a \}}&\mbox{in }\Omega_{\theta_0}.
\end{array}
\right.
\end{eqnarray}
Setting $\lambda = -\sqrt{\frac{{\| \nabla a \|}^2_2+{\| \nabla b \|}^2_2}{2{\| 
\nabla \varphi \|}^2_2}} $ and by arbitrariness of ${\theta_0}$, property 1) is
demonstrated. Now, choosing $\alpha$, $\beta \in C^\infty(\Omega_{\theta_0})$ with
$\alpha=\beta=0$ on $\Gamma_1=\{(r,\theta), \theta= \theta_0\mbox{ or } 
\theta= \theta_0+\frac{2\pi}{m}\}$, it follows from Lemma 2 
$$
\int_{\partial\Omega_{\theta_0}}\frac{ \partial a }{\partial n}\cdot\alpha + \frac{ \partial b }{\partial n}\cdot\beta =0,
$$
that is, $\frac{ \partial a }{\partial n}= \frac{ \partial b }{\partial n}=0$ on $\partial\Omega_{\theta_0}\setminus \Gamma_1$. By arbitrariness of ${\theta_0}$, we establish the property 2).\\
 \noindent
The properties 3) and 4) are just results of 1) , 2) and Lemma 2.\\
\noindent
Now, we choose any vector field $X\in C^\infty(\Omega,\rit^2)$ such that $X\circ A = A\circ X$ and $X\cdot n=0$ on $\partial\Omega$. Let $\sigma_t$ be the flow associated to $X$. Clearly, 
$$
\sigma_t\circ A = A\circ \sigma_t.
$$
Therefore, $\Theta\circ \sigma_t\in F_m$. The proofs of 5) and 6) are the same as the proofs of
Theorem 3.2 (vii) and (iv) in \cite{Ge1}, respectively.
\cqfd

\section{Study of a minimizing sequence}
Through this section we analyze the behaviour of a minimizing sequence in the
spirit of the theory developed in \cite{Ge1} (see also \cite{E}, \cite{PL} or \cite{S2}).
First, we prove a useful fact.\\

\begin{Lem} 
\label{lemma3}
Assume that $\Theta=(a,b) \in F_m$. Then we have
$$
\int_\Omega a = \int_\Omega b =0 .
$$
\cqfd
\end{Lem}
{\em Proof.} By definition of $F_m$, we have
$$
\int_\Omega \Theta = \int_\Omega A^{-1}\circ \Theta\circ A =A^{-1}\int_\Omega  \Theta\circ A = A^{-1}\int_\Omega  \Theta.
$$
We note that $A^{-1}$ is a rotation. This implies 
$$
\int_\Omega \Theta =(0,0).
$$
\cqfd

Now we consider the minimum of energy functional $E$. Set $\ds{G(\Omega) = \inf_{a,b\in H^1\times H^1}E(a,b,\Omega)}$ and $\ds{G_m(\Omega) = \inf_{a,b\in F_m}E(a,b,\Omega)}$. Let $(a_n, b_n, \varphi_n)$ be a minimizing sequence of $E$ on $F_m$, that is, $(a_n, b_n)\in F_m$, $(a_n, b_n,\varphi_n)$ satisfying equation (\ref{eq2}) and
\begin{eqnarray}
\label{eq12bis}
E(a_n, b_n, \Omega) = G_m(\Omega) + o(1).
\end{eqnarray}
Without loss of generality, we can assume that $\|\nabla \varphi_n\|_2=1$. After extracting a subsequence, we may assume that
$$
\begin{array}{l}
a_n \longrightarrow \alpha \mbox{ weakly in }H^1 \mbox{ and strongly in }L^2,\\
b_n \longrightarrow \beta \mbox{ weakly in }H^1 \mbox{ and strongly in }L^2,\\
\varphi_n \longrightarrow \psi \mbox{ weakly in }H^1 \mbox{ and strongly in }L^2.
\end{array}
$$
Obviously, $(\alpha,\beta)\in F_m$. First, we recall a technical lemma.\\

\begin{Lem}
\label{Lemma5}(see \cite{H1}, \cite{Ge1} and also \cite{BC})  We assume that
$\varphi_n$ is a bounded sequence in $H^1_0 \cap L^{\infty}$. Let $a_n
\longrightarrow 0 \mbox{ weakly in }H^1 \mbox{ and strongly in }L^2$. Then for every
$b \in H^1$, we have
\begin{eqnarray}
\lim_{n \rightarrow \infty} \int \varphi_n \{ a_n, b \} =0.
\end{eqnarray}
\cqfd
\end{Lem}

We state the following result, analogous to Theorem 7.1 in \cite{Ge1}.\\

\begin{Lem}
\label{Lemma6}
Under the above assumptions, we have that:\\
(1) if $\psi = 0$, then $\alpha = \beta = 0$;\\
or\\
(2) if $\psi \neq 0$, then $(\alpha, \beta, \psi)$ is a minimum of the energy $E$ on
$F_m$. Moreover, the following holds:
$$
\begin{array}{l}
a_n \longrightarrow \alpha \mbox{  strongly in }H^1,\\
b_n \longrightarrow \beta \mbox{ strongly in }H^1,\\
\varphi_n \longrightarrow \psi \mbox{ strongly in }H^1.
\end{array}
$$
\cqfd
\end{Lem}
{\em Proof.} The proof is the same as the proof of Theorem 7.1 in \cite{Ge1}, but here we
work with equivariant maps.
\cqfd

In the following, we will suppose that $\psi=\alpha=\beta=0$. Denote by $M(\rit^2)$ the space of  non-negative measures on $\rit^2$ with finite mass. Set $\mu_n = \frac{1}{2} (|\nabla a_n|^2 + |\nabla b_n|^2)dx$ and $\nu_n = |\nabla \varphi_n|^2dx$. We consider the extensions of $\mu_n$ and $\nu_n$ to all of $\rit^2$ by valuing 0 in $\rit^2\setminus \Omega$. Then $\{\mu_n\}$ and $\{\nu_n\}$ are bounded in $M(\rit^2)$. Modulo a subsequence, we may assume that
$\mu_n \longrightarrow \mu$, $\nu_n \longrightarrow \nu$ weakly in the sense of measures where $\mu$ and $\nu$ are bounded non-negative measure on $\rit^2$.\\

\begin{Lem}
\label{Lemma7} 
Under assumptions of Lemma 6, if $\psi=\alpha=\beta=0$, then we have
$$
G_m(\Omega)\geq\sqrt{m}G(\Omega).
$$
\cqfd
\end{Lem}

{\em Proof.} Clearly, $\mu(\rit^2\setminus \bar \Omega) = \nu(\rit^2\setminus \bar \Omega) = 0$. Choose $\xi \in C^{\infty}(\rit^2)$. Denote by $ \psi_n$ the unique solution of equation (\ref{eq2}) for $a = \xi a_n$  and $b = \xi b_n$, that is
$$
\left\{
\begin{array}{rll}
- \triangle \psi_n & = \{ \xi a_n, \xi b_n \}, & \mbox{in } \Omega\\
\psi & = 0, & \mbox{on } \partial \Omega.
\end{array}
\right.
$$
A computation using the same arguments as in the proof of Lemma 7.5 in \cite{Ge1} gives
$$
\lim_{n \rightarrow \infty} \| \nabla(\psi_n- \xi^2\varphi_n) \|_2 = 0.
$$
Hence, we obtain
$$
G(\Omega)\| \nabla(\xi^2\varphi_n) \|_2 + o(1) =G(\Omega)\| \nabla\psi_n \|_2 \leq
\frac{1}{2} (
\| \nabla(\xi a_n) \|_2^2 + \| \nabla(\xi b_n) \|_2^2).
$$
Passing to the limit as $n \longrightarrow \infty$, there holds
\begin{eqnarray}
\label{eq15}
G(\Omega)\sqrt{\int \xi^4 d\nu} \leq \int \xi^2 d\mu,\hs\forall \xi \in C^\infty_0(\rit^2). 
\end{eqnarray}
By approximation, therefore,
\begin{eqnarray}
\label{eq16}
G(\Omega)\sqrt{\nu(E)} \leq  \mu(E) \hs (E \subset \rit^2, E \mbox{ Borel}).
\end{eqnarray}
Now since $\mu$ is a finite measure, the set
$$
D \equiv \{ x \in \bar\Omega, \mu(\{ x \}) >0 \}
$$ 
is at most countable. We can therefore write $D=\{x_j\}_{j\in J}$, $\mu_{x_j} = 
\mu(\{x_j\})(j\in J)$ so that
$$
\mu \geq \sum_{j\in J} \mu_{x_j} \delta_{x_j}.
$$
Since (\ref{eq16}) implies $\nu$ is absolutely continuous relative to $ \mu$, we can
write
\begin{eqnarray}
\label{eq17}
\nu(E) = \int_E h d\mu \hs(E \mbox{ Borel}),
\end{eqnarray}
where
\begin{eqnarray}
\label{eq18}
h(x) \equiv \lim_{r\rightarrow 0} \frac{\nu(B(x,r))}{\mu(B(x,r))},
\end{eqnarray}
this limit existing for $\mu$-a.e. $x\in\rit^2$. On the other hand, from (\ref{eq16}), we have
$$
(G(\Omega))^2\frac{\nu(B(x,r))}{\mu(B(x,r))} \leq  \mu(B(x,r)),
$$
provided $\mu(B(x,r))\neq 0$. Thus we infer
\begin{eqnarray}
\label{eq19}
h = 0 \hs \mu\mbox{-a.e. }x\in\rit^2\setminus D.
\end{eqnarray}
Finally, define $\nu_{x_j} \equiv h (x_j) \mu_{x_j}$. Then we get
$$
\nu = \sum_{j\in J} \nu_{x_j} \delta_{x_j},
$$
and
$$
G(\Omega) \sqrt{\nu_{x_j}} \leq\mu_{x_j}.
$$
However, by symmetry of functions in $F_m$, we have $A^i x_j\in D $ for
$i=1,2,...,m-1$ provided $x_j\in D$. Consequently, by suitably relabelling the
$x_j$, we may assume that $x_j\in \Omega_{\theta_0}$ ( where $\Omega_{\theta_0}$ is
defined as in the proof of Lemma 3) for $j\in\{1,...,k\}=J'$ and $k= card(J)m^{-1}$
and
$$
\nu =\sum_{i=0}^{m-1}\sum_{j\in J'} \nu_{x_j} \delta_{A^i x_j} \mbox{ and } \mu \geq \sum_{i=0}^{m-1}\sum_{j\in J'} \mu_{x_j} \delta_{A^i x_j}.
$$
On the other hand, we have $\nu(\bar\Omega)=1$ and $\mu(\bar\Omega)=G_m(\Omega)$. This implies
$$
\begin{array}{ll}
G_m(\Omega)&\ds{ =\mu(\bar\Omega)\geq  m\sum_{j\in J'} \mu_{x_j}\geq m\sum_{j\in J'} G(\Omega)\sqrt{\nu_{x_j}}}\\
&\ds{\geq mG(\Omega)\sqrt{\sum_{j\in J'}\nu_{x_j}} = \sqrt{m}G(\Omega)}.
\end{array}
$$
\cqfd

\section{Proof of the main Theorem}
In view of Lemma \ref{Lemma6}, the result follows if there is no concentration (i.e.
case (1) in Lemma \ref{Lemma6} for minimizing sequences does not occur). By
Lemma \ref{Lemma7}, a sufficient condition for that is to assume
$G_m(\Omega)<\sqrt{m}G(\Omega)$. For this purpose, we set $a(x,y)= x$ and $b(x,y)=y$. It is
obvious to see that $(a,b)\in F_m$ and $E(a,b,\Omega)>0$. For any fixed $r_0>0$, we
can choose some $m\in \nit$ such that 
$$
\sqrt{m}G(\Omega)> E(a,b,\Omega)\geq G_m(\Omega).
$$
Thus, the main Theorem is proved.\\

{\em Acknowlegements.} The authors thank  the referee for his valuable remarks on a
first version of this paper.
\\

\noindent 
Y. Ge ( ge@cmla.ens-cachan.fr )

\noindent 
C.M.L.A., E.N.S de Cachan\\
61, avenue du Pr\'{e}sident Wilson\\
94235 Cachan Cedex, France\\
and\\
D\'epartement de Math\'ematiques\\
Facult\'e de Sciences et  Technologie\\
Universit\'e Paris XII-Val de Marne\\
61, avenue du G\'en\'eral de Gaulle\\
94010 Cr\'eteil Cedex, France

\bigskip

\noindent
F. H\'elein ( helein@cmla.ens-cachan.fr )

\noindent
C.M.L.A., E.N.S de Cachan\\
61, avenue du Pr\'{e}sident Wilson\\
94235 Cachan Cedex, France

\end{document}